\documentclass[a4paper,11pt]{amsart}
\usepackage{etex}
\usepackage{tikz, tikz-cd}

\usepackage{bm}
\usepackage{amssymb}
\usepackage{mathrsfs}
\usepackage{amsmath, amssymb,amsthm,latexsym,amscd,mathrsfs}
\usepackage{indentfirst}
\usepackage{stmaryrd}
\usepackage{graphicx}
\usepackage{subfigure}
\usepackage{extarrows}
\usepackage{amssymb}
\usepackage{amsmath,amssymb,amsthm,latexsym,amscd,mathrsfs}
\usepackage{indentfirst}
\usepackage{multirow}
\usepackage{latexsym}
\usepackage{amsfonts}
\usepackage{color}
\usepackage{pictexwd,dcpic}
\usepackage{graphicx}
\usepackage{psfrag}
\usepackage{hyperref}
\usepackage{comment}
\usepackage{cases}
\usepackage{dutchcal}

\makeatletter
\@namedef{subjclassname@2010}{{\rm2020} Mathematics Subject Classification}
\makeatother

 \newcommand{\ROM}[1]{\mathrm{\uppercase\expandafter{\romannumeral#1}}}
  \theoremstyle{definition}
   \numberwithin{equation}{section}
   \theoremstyle{plain}

 \newtheorem{thm}{Theorem}

 \newtheorem{cor}{Corollary}
  
 \newtheorem{rem}{Remark}

  \makeatletter

  \allowdisplaybreaks
\makeatother
\title[Orthogonal almost complex structure and its Nijenhuis tensor
]{\textbf{Orthogonal almost complex structure and its Nijenhuis tensor}}
\author[Z. Z. Tang]{Zizhou Tang}\address{Chern Institute of Mathematics $\&$ LPMC, Nankai University, Tianjin 300071, P. R. China}
\email{zztang@nankai.edu.cn}
\author[W. J. Yan]{Wenjiao Yan}\address{School of Mathematical Sciences, Laboratory of Mathematics and Complex Systems, Beijing Normal University, Beijing, 100875, P. R. China}
\email{wjyan@bnu.edu.cn}
\thanks {The project is partially supported by the NSFC (No.11931007,12271038), Nankai Zhide Foundation, and the Fundamental Research Funds for the Central Universities (2233300002 )
}

\subjclass[2010] { 32Q60, 53C15, 53C28.}
\keywords{almost Hermitian structure, integrable, Nijenhuis tensor.}

\begin{document}

\maketitle

\begin{abstract}
In this paper, we demonstrate that on an almost Hermitian manifold $(M^{2n}, J, ds^2)$, a 2-form $\varphi=S^*\Phi$, the pulling back of the K\"ahler form $\Phi$ on the twistor bundle over $M^{2n}$, is non-degenerate if the squared norm $|N|^2$ of the Nijenhuis tensor is less than $\frac{64}{5}$ when $n\geq 3$ or less than $16$ when $n=2$. As a corollary, there exists no orthogonal almost complex structure on the standard sphere $(S^6, ds_0^2)$ with $|N|^2<\frac{64}{5}$ everywhere.
\end{abstract}

\section{\textbf{Introduction}}\label{sec1}

Let $(M^{2n}, ds^2)$ be an oriented smooth manifold of dimension $2n$, equipped with a Riemannian metric $ds^2$. An almost complex structure on $M$ is an endomorphism $J$ of the tangent bundle $TM$ such that 
$J^2=-Id$. If the almost complex structure $J$ is an orthogonal transformation with respect to the metric $ds^2$, in other words, $J$ is compatible with the metric: 
\begin{equation}\label{g}
	ds^2(JX, JY)=ds^2(X, Y),~~~\,\,\forall~ X, Y\in TM,
\end{equation} 	
$(M^{2n}, J, ds^2)$ will be called an almost Hermitian manifold. 
On such an $M^{2n}$, one can choose locally a field of $U$-bases $e_1,\cdots, e_{2n}$ satisfying
\begin{equation*}\label{J0}
	J(e_1,e_2,...,e_{2n}) = (e_1,e_2,...,e_{2n})J_0,~~~\mathrm{with}~J_0 =\begin{pmatrix}0&-I_n\\I_n&0\end{pmatrix}.
\end{equation*}
In fact, if we are given an almost complex manifold $(M^{2n}, J)$ first, there always exists a metric $ds^2$ such that $(M^{2n}, J, ds^2)$ is an almost Hermitian manifold. Moreover, equation (\ref{g}) implies that the tensor $\Phi(X, Y)=ds^2(JX, Y)$ is anti-symmetric. $\Phi$ is well known as the fundamental 2-form, or the K\"ahler form on the almost Hermitian manifold $(M^{2n}, J, ds^2)$.

The Newlander-Nirenburg theorem asserts that an almost complex structure $J$ is integrable if and only if the Nijenhuis tensor
\begin{equation}\label{N}
	N(X, Y)=[JX, JY]-[X, Y]-J[JX, Y]-J[X, JY],~~~\, X, Y\in\mathfrak{X}(M)
\end{equation}
vanishes everywhere. In this circumstance, the almost complex manifold $(M^{2n}, J)$ would admit local complex coordinates, transforming it into a complex manifold, and the almost Hermitian manifold $(M^{2n}, J, ds^2)$ is referred to as a Hermitian manifold.

After a series of significant contributions by H. Hopf \cite{Hop48}, C. Ehresmann \cite{Ehr49}, W.-T. Wu \cite{Wu52}, A. Borel and J.-P. Serre \cite{BS53}, and others, the renowned Hopf problem has evolved into the problem of determining the existence of a complex structure on $S^6$. Many esteemed mathematicians, including S. S. Chern, M. F. Atiyah, S. T. Yau, etc., have devoted considerable attention to this problem and made noteworthy contributions.
Irrespective of the existence of a complex structure on $S^6$, the authors  \cite{TY22} found an
$8$-dimensional closed manifold $N^8$ such that a complex structure exists on $S^6\times N^8$.

In 1953, A. Blanchard \cite{Bla53} showed that there are no orthogonal almost complex structures $J$ on the standard sphere $(S^6, ds^2_0)$ which is integrable (equivalently, the Nijenhuis tensor $N\equiv 0$). This finding was reaffirmed by LeBrun in 1987  \cite{LeB87}, and was considered as a milestone in the context of the Hopf problem. As the primary outcome of this paper, we generalize this result in the following manner.

Our strategy was developed based on the twistor bundle theory (see, for example, \cite{PT04}), details will be explained in Section \ref{2-form}.
Given an oriented Riemannian manifold $(M^{2n}, ds^2)$, we explore the positive principal frame bundle $SO(2n)\hookrightarrow P_+(M, ds^2)\xrightarrow{p} M$, 
and the twistor bundle $\Gamma_n=SO(2n)/U(n)\hookrightarrow J_+(M, ds^2)\xrightarrow{\pi}M$. 
There exist natural projections in this context
\begin{equation*}\label{pb}
	\begin{array}{lll}
		P_+(M, ds^2)&\xlongrightarrow{\xi}& J_+(M, ds^2)\\
		\qquad\quad p\searrow &&\swarrow \pi\\
		&M&
	\end{array}
\end{equation*}
with $\pi\circ\xi=p$, where $\xi: P_+(M, d^2)\rightarrow J_+(M, ds^2)$ is a $U(n)$-principal bundle. Clearly, $J_+(M, ds^2)$ is the associated fiber bundle of the principal bundle $p$.                                                                                                                                                                                                                                                                                                                                                                                                                                                                                                                                                                                                                                                                                                                                                                                                                                                                                                                                                                                                                                                                                                                                                                                                                                                                                                                                                                                                                                                                                                                                                                                                                                   
For each $e\in P_+(M, ds^2)$, a linear isomorphism can be defined as follows:
\begin{eqnarray}\label{lambda1}
	T_eP_+(M, ds^2)=V_e\oplus H_e& \longrightarrow& \mathfrak{so}(2n)\oplus \mathbb{R}^{2n}\\
	\widetilde{X}\quad\,\,\,&\mapsto& (~\Omega(\widetilde{X}), ~~\Theta(\widetilde{X})~)\nonumber
\end{eqnarray}
where $\Omega=(\Omega_{AB})_{2n\times 2n}$ are the $\mathfrak{so}(2n)$-valued Levi-Civita connection $1$-forms on the principal bundle $P_+(M, ds^2)$, and $\Theta=(\Theta_1,\cdots,\Theta_{2n})$ are the $\mathbb{R}^{2n}$-valued canonical 1-forms on $P_+(M, ds^2)$. Decompose 
the Lie algebra $\mathfrak{so}(2n)$ as $\mathfrak{so}(2n)=\mathfrak{u}(n)\oplus \Sigma$, where $\mathfrak{u}(n)$ is the Lie algebra of $U(n)$. 
Inheriting from the isomorphism in (\ref{lambda1}), a well-defined isomorphism  can be given:
\begin{eqnarray*}
	\lambda(e): ~T_{\xi(e)}J_+(M, ds^2)&\longrightarrow& \Sigma\oplus \mathbb{R}^{2n}\\
	\bar{X}\qquad&\mapsto&(~\frac{1}{2}(\Omega+J_0\Omega J_0)(\widetilde{X}),~\Theta(\widetilde{X})~)\nonumber
\end{eqnarray*}
where $\widetilde{X}$ is an element in $T_eP_+(M, ds^2)$ such that $\xi_*(\widetilde{X})=\bar{X}$. 
By pulling back the canonical complex structure and the canonical metric on $\Sigma\oplus \mathbb{R}^{2n}$, one derives an almost complex structure $J_1$ on $J_+(M, ds^2)$, and a metric $ds^2_{J_1}$ that is compatible with $J_1$. Consequently, one obtains an almost Hermitian manifold $(J_+(M, ds^2), J_1, ds_{J_1}^2)$, where the K\"ahler form is $\Phi(\bar{X}, \bar{Y})=ds_{J_1}^2(J_1\bar{X}, \bar{Y})$.

Notice that every almost complex structure $J$ on $M$ corresponds to a cross section of the twistor bundle $\pi: J_+(M, ds^2)\rightarrow M$, denoted by $S$.  Pulling back the K\"ahler form $\Phi$,  
we obtain a 2-form $\varphi$ on $M$:
\begin{equation}\label{varphi}
	\varphi(X, Y):=S^*\Phi(X, Y)=\Phi(S_*X, S_*Y), \,\,\, X, Y\in TM
\end{equation}

As usual, the squared norm $|N|^2$ of the Nijenhuis tensor is defined by $|N|^2:=\sum_{A,B=1}^{2n}|N(e_A, e_B)|^2$. 
As the main result of this paper, we establish
\begin{thm}\label{thm}
	On an almost Hermitian manifold $(M^{2n}, J, ds^2)$, the 2-form $\varphi$ in (\ref{varphi}) is non-degenerate if $|N|^2$ is less than $c_0$, where
$c_0=\frac{64}{5}$ when $n\geq 3$ and $c_0=16$ when $n=2$.
\end{thm}

As a consequence, we generalize the result of Blanchard and LeBrun:
\begin{cor}\label{cor}
	On the standard sphere $(S^6, ds_0^2)$, there is no orthogonal almost complex structure $J$ 
	with $|N|^2<\frac{64}{5}$ everywhere.
\end{cor}

\begin{rem}
In \cite{LeB87}, an almost complex structure on the standard $S^6$ is shown to define a map from $S^6$ to the Grassmannian
$G_3(\mathbb{C}^7)$ which is a K\"ahler manifold. If the almost complex structure is integrable, the map would embed $S^6$ as
a complex manifold into $G_3(\mathbb{C}^7)$, giving $S^6$ a K\"ahler structure. This is impossible since $H^2(S^6)=0$.
However, LeBrun's method depends on the K\"ahler structure very much, and is too rigid to get our Corollary \ref{cor}.
\end{rem}

The organization of this paper is as follows: In Section \ref{2-form}, we introduce the details about $P_+(M, ds^2)$ and $J_+(M, ds^2)$, and give the precise expression of $\varphi$; in Section \ref{3}, we provide the proof of Theorem  \ref{thm} and Corollary \ref{cor}.


\section{The 2-form $\varphi$}\label{2-form}
We first revisit some foundational concepts from \cite{PT04}. As usual, interpret $U(n)$ as a closed subgroup of $SO(2n)$: $U(n)=\{A\in SO(2n)~|~AJ_0=J_0A\}$. 
On an oriented Riemannian manifold $(M^{2n}, ds^2)$, $(T_xM, ds_x^2)$ is an oriented real vector space for any $x\in M$. Let 
$$ P_+(T_xM, ds_x^2)=\{e=(e_1, \cdots, e_{2n})~|~e_j\in T_xM, ds_x^2(e_i, e_j)=\delta_{ij}, e>0 \}$$ 
be the space of orthonormal frames with positive orientation, and 
\begin{eqnarray*}
J_+(T_xM, ds_x^2)&=&\{J\in \mathrm{Hom}(T_xM, T_xM)~|~J^2=-Id, \\
&&J \text{ is compatible with } \text{the metric and the orientation} \}
\end{eqnarray*}
be the set of ``positive" Hermitian structures on $(T_xM, ds_x^2)$. For any $e\in P_+(T_xM, ds_x^2)$, there exists a natural onto projection $\xi: P_+(T_xM, ds_x^2)\rightarrow J_+(T_xM, ds_x^2)$ defined by $$\xi(e)(e_1,\cdots,e_{2n})=(e_1, \cdots, e_{2n})J_0.$$ It is easy to see that $\xi(e)=\xi(e')$ if and only if $e=e'U$ for some $U\in U(n)$.

Set
$$	P_+(M^{2n}, ds^2)=\underset{x\in M}{\cup} P_+(T_xM, ds_x^2),$$
and 
$$	J_+(M^{2n}, ds^2)=\underset{x\in M}{\cup} J_+(T_xM, ds_x^2).$$ 
It is well known that $p: P_+(M, ds^2)\rightarrow M$ is a ``positive" principal frame bundle with the structure group $SO(2n)$ acting on $P_+(M, ds^2)$ on the right, and $\xi: P_+(M, ds^2)\rightarrow J_+(M, ds^2)$ is a $U(n)$-principal bundle.
Moreover, there are natural projections
\begin{equation}\label{pb}
	\begin{array}{lll}
		P_+(M, ds^2)&\xlongrightarrow{\xi}& J_+(M, ds^2)\\
		\qquad\quad p\searrow &&\swarrow \pi\\
		&M&
	\end{array}
\end{equation}
so that $\pi\circ\xi=p$.
It is known that $J_+(M, ds^2)$ is the associated fiber bundle of the principal bundle $p$ with fiber $\Gamma_n=SO(2n)/U(n)$. More precisely, 
based on the right action of $SO(2n)$ on $P_+(T_xM, ds_x^2)$ and the left action of $SO(2n)$ on $\Gamma_n$, a left action of $SO(2n)$ on the product $P_+(T_xM, ds_x^2)\times \Gamma_n$ can be defined as $A\cdot(e, BU(n))=(eA^{-1}, ABU(n))$ for $A, B\in SO(2n)$, $e\in P_+(T_xM, ds_x^2)$. Consequently, one obtains a one-to-one correspondence between the quotient space 
$P_+(T_xM, ds_x^2)\times_{SO(2n)} \Gamma_n$ and $J_+(T_xM, ds_x^2)$ through the map $P_+(T_xM, ds_x^2)\times \Gamma_n\rightarrow J_+(T_xM, ds_x^2)$, ~$(e, BU(n))\mapsto \xi(eB)$.


As pointed out in the introduction, we will introduce an almost Hermitian manifold $(J_+(M, ds^2), J_1, ds^2_{J_1})$. To do this, we first investigate the total space $P_+(M, ds^2)$ of the $U(n)$-principal bundle over $J_+(M, ds^2)$.
For $e\in P_+(M, ds^2)$, denote by $V_e=dp_e$ the vertical subspace of $T_eP_+(M, ds^2)$. The connection in the principal bundle $P_+(M, ds^2)$ provides an assignment of the horizontal subspace $H_e$ such that $T_eP_+(M, ds^2)=V_e\oplus H_e$. For any $\widetilde{X}\in T_eP_+(M, ds^2)$, define $\Omega(\widetilde{X})$ to be the unique $l\in\mathfrak{so}(2n)$ such that $l^*(e)$ is the vertical component of $\widetilde{X}$, where $l^*$ is the fundamental vector field induced by the action of $\mathrm{exp}(tl)$ on $P_+(M, ds^2)$. As a result, we obtain the  $\mathfrak{so}(2n)$-valued Levi-Civita connection 
$1$-forms $\Omega=(\Omega_{AB})_{2n\times 2n}$ on the principal bundle $P_+(M, ds^2)$. On the other hand, let $\Theta=(\Theta_1,\cdots,\Theta_{2n})$ be the canonical 1-forms on $P_+(M, ds^2)$ defined by $\Theta_i(\widetilde{X})=ds^2(p_*\widetilde{X}, e_i)$ $(i=1,\cdots,2n)$. In this way, we establish a linear isomorphism for each $e\in P_+(M, ds^2)$:
\begin{eqnarray}\label{lambda}
T_eP_+(M, ds^2)=V_e\oplus H_e& \longrightarrow& \mathfrak{so}(2n)\oplus \mathbb{R}^{2n}\\
\widetilde{X}\quad\,\,\,&\mapsto& (~\Omega(\widetilde{X}), \Theta(\widetilde{X})~).\nonumber
\end{eqnarray}

Next, we decompose the Lie algebra $\mathfrak{so}(2n)$ as $\mathfrak{so}(2n)=\mathfrak{u}(n)\oplus \Sigma$,
where $\mathfrak{u}(n)=\{\phi\in\mathfrak{so}(2n)~|~\phi J_0=J_0\phi\}$ is the Lie algebra of $U(n)$, and $\Sigma$ can be expressed as $\Sigma=\{\psi\in\mathfrak{so}(2n)~|~\psi J_0=-J_0\psi\}$. Consequently, we have the following decomposition for any $\Omega\in\mathfrak{so}(2n)$:
$$\Omega=\frac{1}{2}(\Omega-J_0\Omega J_0)+\frac{1}{2}(\Omega+J_0\Omega J_0)\in \mathfrak{u}(n)\oplus \Sigma.$$
Inheriting from the isomorphism in (\ref{lambda}), we can define the following map for each $e\in P_+(M, ds^2)$: 
\begin{eqnarray}\label{lambda e}
	\lambda(e): ~T_{\xi(e)}J_+(M, ds^2)&\longrightarrow& \Sigma\oplus \mathbb{R}^{2n}\\
	\bar{X}\qquad&\mapsto&(~\frac{1}{2}(\Omega+J_0\Omega J_0)(\widetilde{X}),~\Theta(\widetilde{X})~)\nonumber
\end{eqnarray}
where $\widetilde{X}$ is an element in $T_eP_+(M, ds^2)$ such that $\xi_*(\widetilde{X})=\bar{X}$. It is straightforward to verify that $\lambda(e)$ is a well-defined isomorphism.

On the direct sum $\Sigma\oplus \mathbb{R}^{2n}$, there is a right action of $U(n)$ on $ \Sigma\oplus \mathbb{R}^{2n}$ that is induced by the right action of $U(n)$ on $\Sigma$ by the adjoint and the right action of $U(n)$ on $\mathbb{R}^{2n}$ by the inclusion $U(n)\hookrightarrow SO(2n)$. 
It is easy to see that $\lambda(e)$ is compatible with the right action of $U(n)$ on $\Sigma\oplus \mathbb{R}^{2n}$:
$$\lambda(eU)(\bar{X})=U\cdot \lambda(e)(\bar{X}).$$ 
This guarentees that we can construct the almost complex structure and metric on $J_+(M, ds^2)$ by virtue of $\lambda(e)$. There exists a canonical complex structure expressed as:
\begin{equation}\label{canonical complex str}
	(\psi, V)\longmapsto (J_0\psi, -VJ_0).
\end{equation}
The right action of $U(n)$ on it is holomorphic. Pulling back this complex structure by $\lambda(e)$ in (\ref{lambda e}) to $T_{\xi(e)}J_+(M, ds^2)$, we obtain a canonical almost complex structure $J_1$ on $J_+(M, ds^2)$.
Interested readers can refer to Eells-Salamon \cite{ES85}. 
Furthermore, pulling bcak the canonical metric of $\Sigma\oplus \mathbb{R}^{2n}$ to $( J_+(M, ds^2), J_1)$ and denoting it by $ds^2_{J_1}$, one obtains an almost Hermitian manifold $(J_+(M, ds^2), J_1, ds_{J_1}^2)$, where we have the K\"ahler form $\Phi(\bar{X}, \bar{Y})=ds_{J_1}^2(J_1\bar{X}, \bar{Y})$.

Given a positive orthogonal almost complex structure $J\in J_+(M^{2n}, ds^2)$ with corresponding cross section $S$, by pulling back the K\"ahler form $\Phi$ of $(J_+(M, ds^2), J_1, ds_{J_1}^2)$ to $M$,  
we obtain a 2-form $\varphi$ on $M$ as expressed in (\ref{varphi}):
\begin{equation*}
	\varphi(X, Y):=S^*\Phi(X, Y)=ds_{J_1}^2(J_1S_*X, S_*Y),\,\,\, X, Y\in TM.
\end{equation*}

At last, let's interpret the 2-form $\varphi$ in terms of $\Omega$ and $\Theta$. Choose a local cross section corresponding to $\xi$ in $(\ref{pb})$ by $\mu$, such that $\mu\circ S=e$, 
and denote $\omega=e^*\Omega$ and $\theta=e^*\Theta$. Therefore, $\theta\theta^T=ds^2$, and $\omega$ is just the Levi-Civita connection form on $M$. 
They satisfy the structure equations:
$$\begin{cases} 
	\rm{d}\theta=\theta\wedge\omega.\\
	\rm{d}\omega=\omega\wedge\omega-R.
\end{cases}
$$
Then it follows from the correspondence in $(\ref{lambda})$ that
\begin{eqnarray*}
S_*X &\longleftrightarrow& (~\frac{1}{2}(\Omega+J_0\Omega J_0)(\mu_*S_*X),~~\Theta(\mu_*S_*X))\\
&&=(~\frac{1}{2}e^*(\Omega+J_0\Omega J_0)(X),~~e^*\Theta(X))\\
&&=(~\frac{1}{2}(\omega+J_0\omega J_0)(X),~~\theta(X)),
\end{eqnarray*}
and 
\begin{eqnarray*}
	J_1S_*X &\longleftrightarrow&(~\frac{1}{2}J_0(\omega+J_0\omega J_0)(X),~~-\theta(X)J_0~)\\
	&&=(~\frac{1}{2}(J_0\omega-\omega J_0)(X),~~-\theta(X)J_0~).
\end{eqnarray*}
Therefore, 
\begin{eqnarray}\label{phi}
	\varphi(X, Y)&=&ds_{J_1}^2(J_1S_*X, S_*Y)\nonumber\\
	&=&\frac{1}{4} ds_{\mathfrak{so}(2n)}^2\big((J_0\omega-\omega J_0)(X), (\omega+J_0\omega J_0)(Y)\big)-ds_{{\mathbb{R}^{2n}}}^2\big(\theta(X)J_0, \theta(Y) \big).
\end{eqnarray}

Now, on the almost Hermitian manifold $(M^{2n}, J, ds^2)$, using indices $ i, j, k=1,\cdots,n$, $A, B, C=1,\cdots, 2n$, 
we choose locally an orthonormal frame field $\{e_1, e_2, \cdots, e_{2n}\}$ satisfying
$$e_{n+1}=Je_1,~~ e_{n+2}=J e_2,~~\cdots, ~~e_{2n}=Je_n.$$
For simplicity, let's denote 
\begin{eqnarray}\label{ab}
	\alpha_{ij}:=\omega_{i, j+n}+\omega_{i+n, j},&&\alpha_{ij}^k:=\alpha_{ij}(e_k),\\
	\beta_{ij}:=\omega_{i+n, j+n}-\omega_{ij},~&&\beta_{ij}^k:=\beta_{ij}(e_k)\nonumber.
\end{eqnarray}
Then 
$$J_0\omega-\omega J_0=\begin{pmatrix} -\omega_{i, j+n}-\omega_{i+n, j}&-\omega_{i+n, j+n}+\omega_{ij}\\ -\omega_{i+n, j+n}+\omega_{ij}&\omega_{i, j+n}+\omega_{i+n, j}\end{pmatrix} =\begin{pmatrix} -\alpha_{ij}&-\beta_{ij}\\ -\beta_{ij}&\alpha_{ij}\end{pmatrix} ,$$
$$\omega+J_0\omega J_0=\begin{pmatrix} \omega_{i, j}-\omega_{i+n, j+n}&\omega_{i, j+n}+\omega_{i+n,j}\\ \omega_{i, j+n}+\omega_{i+n,j}&\omega_{i+n, j+n}-\omega_{ij}\end{pmatrix} =\begin{pmatrix} -\beta_{ij}&\alpha_{ij}\\ \alpha_{ij}&\beta_{ij}\end{pmatrix}.$$
Therefore, we derive the precise expression of $\varphi$:
\begin{equation}\label{varphi2}
		\varphi(X, Y)=\frac{1}{2} \sum_{i,j=1}^n \alpha_{ij}\wedge\beta_{ij} (X,Y)+\sum_{i=1}^n\theta_i\wedge\theta_{i+n}(X, Y),
\end{equation}
where we take $ds_{\mathfrak{so}(2n)}^2(P,Q)=-\mathrm{tr}(PQ)$ for convenience.

\vspace{8mm}

\section{Proof of Theorem \ref{thm} and Corollary \ref{cor}}\label{3}

\noindent
\emph{Proof of Theorem \ref{thm}}.

We will establish the non-degeneracy of $\varphi$ by showing that $\varphi(X, JX)>0$ for any $X\in TM$ with $X\neq 0$.

Express $X\in TM$ as 
$X=\underset{k}\sum \big(x_ke_k+x_{k+n}e_{k+n}\big)$. It is evident that 
$$\sum_{i=1}^n\theta_i\wedge\theta_{i+n}(X, JX)=\sum_i\theta_i(X)\theta_{i+n}(JX)-\theta_i(JX)\theta_{i+n}(X)=|X|^2.$$
For simplicity, let's denote 
\begin{eqnarray}\label{Cd}
	C_{ijk}:=\alpha_{jk}^{i+n}+\beta_{jk}^i,~~~\,\,\,~~~&& C'_{ijk}:=\alpha_{jk}^i-\beta_{jk}^{i+n}\\
	d_{ijk}:=C_{ijk}-C_{jik},~~~\,\,\,~~~&&d'_{ijk}:=C'_{ijk}-C'_{jik}\label{d}\nonumber
\end{eqnarray}
Then we have
\begin{eqnarray*}
&&	\alpha_{ij}\wedge\beta_{ij}(X, JX)\\
	&=& \alpha_{ij}(X)\beta_{ij}(JX)-\alpha_{ij}(JX)\beta_{ij}(X)\\
	&=&\big(\alpha_{ij}(X)-\beta_{ij}(JX)\big)\beta_{ij}(JX)-\big(\alpha_{ij}(JX)+\beta_{ij}(X)\big)\beta_{ij}(X)+\beta_{ij}(JX)^2+\beta_{ij}(X)^2.
\end{eqnarray*}
Notice that 
$$\alpha_{ij}(X)-\beta_{ij}(JX)=\sum_k x_k(\alpha_{ij}^k-\beta_{ij}^{k+n})+x_{k+n}(\alpha_{ij}^{k+n}+\beta_{ij}^k)=\sum_k x_kC'_{kij}+x_{k+n}C_{kij},$$
we have
$$\big(\alpha_{ij}(X)-\beta_{ij}(JX)\big)\beta_{ij}(JX)\geq -\sqrt{\sum_k(C_{kij}^2+C'^2_{kij})}~|X|~|\beta_{ij}(JX)|.$$
Similarly,
$$-\big(\alpha_{ij}(JX)+\beta_{ij}(X)\big)\beta_{ij}(X)\geq -\sqrt{\sum_k(C_{kij}^2+C'^2_{kij})}~|X|~|\beta_{ij}(X)|.$$
Denoting 
\begin{equation}\label{A}
	A_{ij}:=\sqrt{\sum_k(C_{kij}^2+C'^2_{kij})},
\end{equation} we derive that
\begin{eqnarray*}
	\alpha_{ij}\wedge\beta_{ij}(X, JX)&\geq&\beta_{ij}(JX)^2-A_{ij}~|X|~|\beta_{ij}(JX)|+\beta_{ij}(X)^2-A_{ij}~|X|~|\beta_{ij}(X)|\\
	&=&\big(\beta_{ij}(JX)-\frac{1}{2}A_{ij}|X|\big)^2+\big(\beta_{ij}(X)-\frac{1}{2}A_{ij}|X|\big)^2-\frac{1}{2}A_{ij}^2|X|^2\\
	&\geq& -\frac{1}{2}A_{ij}^2|X|^2.
\end{eqnarray*}
Therefore,
\begin{eqnarray}\label{XJX}
		\varphi(X, JX)&=&\frac{1}{2} \sum_{i,j=1}^n \alpha_{ij}\wedge\beta_{ij} (X,JX)+\sum_{i=1}^n\theta_i\wedge\theta_{i+n}(X, JX)\\
		&\geq&|X|^2\big(1-\frac{1}{4}\sum_{i,j}A_{ij}^2\big).\nonumber
\end{eqnarray}

Next, we will estimate $\sum_{i,j}A_{ij}^2$ via the squared norm of the Nijenhuis Tensor $N$.
Using notations in (\ref{ab}) and (\ref{Cd}), $N(e_i, e_j)$
 can be expressed as
 \begin{eqnarray*}
	N(e_i, e_j)&=&[Je_i, Je_j]-[e_i, e_j]-J[Je_i, e_j]-J[e_i, Je_j]\\
	&=&(\nabla_{Je_i}Je_j-J\nabla_{Je_i}e_j)-(\nabla_{Je_j}Je_i-J\nabla_{Je_j}e_i)-(\nabla_{e_i}e_j+J\nabla_{e_i}Je_j)\\
	&&+(\nabla_{e_j}e_i+J\nabla_{e_j}Je_i)\\
	&=&\sum_A \Big(\omega_{j+n, A}(Je_i)e_A-J\omega_{jA}(Je_i)e_A\Big)-\sum_A \Big(\omega_{i+n, A}(Je_j)e_A-J\omega_{iA}(Je_j)e_A\Big)\\
	&&-\sum_A \Big(\omega_{j, A}(e_i)e_A+J\omega_{j+n,A}(e_i)e_A\Big)+\sum_A \Big(\omega_{i, A}(e_j)e_A+J\omega_{i+n,A}(e_j)e_A\Big)\\
&=& \sum_k  \Big\{ \big(\alpha_{jk}^{i+n}+\beta_{jk}^i-\alpha_{ik}^{j+n}-\beta_{ik}^j\big)e_k
+\big(\beta_{jk}^{i+n}-\alpha_{jk}^i-\beta_{ik}^{j+n}+\alpha_{ik}^j\big)e_{k+n}\Big\}\\
&=& \sum_k \Big\{ \big( C_{ijk}-C_{jik}\big)e_k -\big( C'_{ijk}-C'_{jik}\big)e_{k+n}\Big\}\\
&=& \sum_k \big(d_{ijk}~e_k-d'_{ijk}~e_{k+n}\big).
\end{eqnarray*}
Thus
\begin{equation*}\label{N2}
	|N(e_i, e_j)|^2=\sum_k \big(d_{ijk}^2+d'^2_{ijk}\big).
\end{equation*}
Recalling the properties that
$N(Y, X)=-N(X, Y)$, $N(JX, Y)=-JN(X, Y)=N(X, JY)$, 
we derive the squared norm of the Nijenhuis tensor as follows:
\begin{equation}\label{N3}
|N|^2:=\sum_{A,B}|N(e_A, e_B)|^2=4\sum_{i,j}|N(e_i, e_j)|^2=4\sum_{i,j,k}\big(d_{ijk}^2+d'^2_{ijk}\big) .
\end{equation}

At last, we divide the investigation into two cases.

\noindent
\textbf{Case 1: $n\geq 3$.}
From the notations in (\ref{ab}) and (\ref{Cd}), it follows that
\begin{eqnarray*}
	C_{ijk}&=&-C_{ikj}=-C_{kij}-d_{ikj}=C_{kji}-d_{ikj}\\
	&=&C_{jki}+d_{kji}-d_{ikj}=-C_{jik}+d_{kji}-d_{ikj}\\
	&=&-C_{ijk}-d_{jik}+d_{kji}-d_{ikj},
\end{eqnarray*}
that is,
\begin{equation}\label{C1}
	2C_{ijk}=-d_{jik}+d_{kji}-d_{ikj}=d_{ijk}-d_{jki}+d_{kij}.
\end{equation}
Similarly, we have
\begin{equation}\label{C2}
2C_{jki}=d_{jki}-d_{kij}+d_{ijk},\quad 2C_{kij}=d_{kij}-d_{ijk}+d_{jki}.
\end{equation}
Squaring both sides of (\ref{C1}) and (\ref{C2}) and summing them together, we acquire
\begin{eqnarray*}
4(C_{ijk}^2+C_{jki}^2+C_{kij}^2)&=&3(d_{ijk}^2+d_{jki}^2+d_{kij}^2)-2d_{ijk}d_{jki}-2d_{ijk}d_{kij}-2d_{jki}d_{kij}\\
&\leq&5(d_{ijk}^2+d_{jki}^2+d_{kij}^2).
\end{eqnarray*}
Consequently, $$\sum_{i,j,k}C_{ijk}^2\leq\frac{5}{4}\sum_{i,j,k}d^2_{ijk}.$$
Analogously,
$$\sum_{i,j,k}C'^2_{ijk}\leq\frac{5}{4}\sum_{i,j,k}d'^2_{ijk}.$$
Therefore, if $|N|^2<\frac{64}{5}$, then it follows from (\ref{A}), (\ref{XJX}) and (\ref{N3}) that
\begin{eqnarray*}
\varphi(X, JX)&\geq&|X|^2\big(1-\frac{1}{4}\sum_{i,j,k}(C^2_{ijk}+C'^2_{ijk})\big)
\geq|X|^2\big(1-\frac{5}{16}\sum_{i,j,k}(d^2_{ijk}+d'^2_{ijk})\big)\\
&=&|X|^2\big(1-\frac{5}{64}|N|^2\big)>0
\end{eqnarray*}
for any non-vanishing $X\in TM$.
\vspace{5mm}

\noindent
\textbf{Case 2: $n=2$.}
In this case, one can calculate directly that $$\sum_{i,j,k}C_{ijk}^2=2(d^2_{121}+d^2_{212})=\sum_{i,j,k}d_{ijk}^2,$$
and $$\sum_{i,j,k}C'^2_{ijk}=2(d'^2_{121}+d'^2_{212})=\sum_{i,j,k}d'^2_{ijk}.$$
Therefore, if $|N|^2<16$, then it follows from (\ref{A}), (\ref{XJX}) and (\ref{N3}) that
\begin{eqnarray*}
	\varphi(X, JX)&\geq&|X|^2\Big(1-\frac{1}{4}\sum_{i,j,k}(C^2_{ijk}+C'^2_{ijk})\Big)
	= |X|^2\Big(1-\frac{1}{4}\sum_{i,j,k}(d^2_{ijk}+d'^2_{ijk})\Big)\\
	&=&|X|^2\big(1-\frac{1}{16}|N|^2\big)>0
\end{eqnarray*}
for any non-vanishing $X\in TM$.

\hfill
$\Box$

\noindent
\emph{Proof of Corollary \ref{cor}.}

From the Cartan structure equation, we have
\begin{eqnarray}
	\sum_id\omega_{i, i+n}&=&\sum_{i, A}\omega_{iA}\wedge\omega_{A,i+n}-\sum_iR_{i,i+n}\\
	&=&\sum_{i,j}\Big(\omega_{ij}\wedge\omega_{j,i+n}+\omega_{i,j+n}\wedge\omega_{j+n, i+n}\Big)-\sum_iR_{i,i+n}.\nonumber
\end{eqnarray}
It is stated in \cite{Tan06} that $\sum_id\omega_{i, i+n}$ is a global 2-form on the almost Hermitian manifold
$(M^{2n}, J, ds^2)$, and $-\frac{1}{2\pi}\sum_id\omega_{i, i+n}$ represents the first Chern class.
Through a direct calculation, one can verify that:
$$\sum_{i,j}\Big(\omega_{ij}\wedge\omega_{j,i+n}+\omega_{i,j+n}\wedge\omega_{j+n, i+n}\Big)=-\frac{1}{2}\sum_{i,j}\alpha_{ij}\wedge\beta_{ij}.$$
Therefore, located on the standard sphere $(S^6, ds_0^2)$, we have
\begin{eqnarray*}
\sum_id\omega_{i, i+n}&=&-\frac{1}{2}\sum_{i,j}\alpha_{ij}\wedge\beta_{ij}-\sum_i\theta_i\wedge\theta_{i+n}+\sum_i\theta_i\wedge\theta_{i+n}-\sum_iR_{i,i+n}\\
&=&-\varphi,
\end{eqnarray*}
which implies that $\varphi$ is a closed 2-form on $(S^6, ds_0^2)$. Since $H^2(S^6; \mathbb{R})\cong0$, there exists $\psi\in\wedge^1(S^6)$ such that $\varphi=d\psi$. Consequently, by the Stokes formula, we have $\int_{S^6}\varphi^3=0$.

However, when $|N|^2<\frac{64}{5}$, it follows from Theorem \ref{thm} that $\varphi$ is non-degenerate. Hence
$\int_{S^6}\varphi^3\neq 0$. There comes a contradiction.

\hfill$\Box$


\end{document}